\documentclass[review,onefignum,onetabnum]{siamart190516}

\usepackage{lipsum}
\usepackage{amsfonts}
\usepackage{graphicx}
\usepackage{epstopdf}
\usepackage{algorithmic}

\usepackage{tcolorbox}
\usepackage{color}

\ifpdf
  \DeclareGraphicsExtensions{.eps,.pdf,.png,.jpg}
\else
  \DeclareGraphicsExtensions{.eps}
\fi

\newsiamremark{hypothesis}{Hypothesis}
\crefname{hypothesis}{Hypothesis}{Hypotheses}
\newsiamthm{claim}{Claim}

\title{Efficient Shallow Ritz Method \\ [1mm] for 1D diffusion-reaction problems\thanks{This work was supported in part by the National Science Foundation under grant DMS-2110571 and by the Department of Energy under NO.B665416. This work was performed under the auspices of the U.S. Department of Energy by Lawrence Livermore National Laboratory under Contract DE-AC52-07NA27344 (LLNL-JRNL-865920).}}

\author{Zhiqiang Cai\thanks{Department of Mathematics, Purdue University, West Lafayette, IN
    (\email{caiz@purdue.edu},  \email{adoktoro@purdue.edu}, \email{herre125@purdue.edu}).}
\and Anastassia Doktorova\footnotemark[2]
\and Robert D. Falgout \thanks{Lawrence Livermore National Laboratory, Livermore, CA (\email{rfalgout@llnl.gov})}
\and César Herrera\footnotemark[2]}

\usepackage{amsopn}

\ifpdf
\hypersetup{
  pdftitle={fastNNSolver},
  pdfauthor={}
}
\fi

\usepackage{amsmath,rotating}

\usepackage[utf8]{inputenc}
\usepackage{bm}
\usepackage{bbm}
\usepackage{mathtools}
\usepackage{indentfirst}
\usepackage{subfigure}
\usepackage{tcolorbox}
\usepackage{color}
\usepackage[export]{adjustbox}
\usepackage{multirow}
\usepackage{scalerel}
\usepackage{tikz}
\usepackage{geometry}
\usepackage{multicol}
\usetikzlibrary{decorations.pathreplacing}
\usepackage{graphicx}
\usepackage{subcaption}
\usepackage{amsmath,amssymb}
\usepackage{nicematrix}
\usepackage{float}

\usepackage{algorithmic}
\Crefname{ALC@unique}{Line}{Lines}

\newcommand{\R}{\mathbb{R}}

\usepackage{graphicx}
\usepackage{diagbox}
\usepackage{pifont}
\definecolor{RED}{rgb}{1, 0, 0} 
\usepackage{soul}
\newcommand{\vertiii}[1]{{\left\vert\kern-0.25ex\left\vert\kern-0.25ex\left\vert #1 
    \right\vert\kern-0.25ex\right\vert\kern-0.25ex\right\vert}}

\newcommand{\bphi}{\mbox{\boldmath${\varphi}$}}
\newcommand{\bpsi}{\mbox{\boldmath${\psi}$}}

\newcommand{\bxi}{\mbox{\boldmath$\xi$}}

\newtheorem{rem}{\bf Remark}[section]

\usepackage{enumitem}
\setlist[itemize]{left=16pt} 

\def\br{{\bf r}}
\def\bb{{\bf b}}
\def\bB{{\bf B}}
\def\ba{{\bf a}}
\def\bc{{\bf c}}
\def\bd{{\bf d}}
\def\bD{{\bf D}}
\def\bH{{\bf H}}
\def\bh{{\bf h}}
\def\bff{{\bf f}}

\def\bv{{\bf v}}

\def\bw{{\bf w}}

\def\bs{{\bf s}}
\def\bg{{\bf g}}
\def\bq{{\bf q}}

\def\cM{{\cal M}}

\begin{document}


%

\maketitle
\begin{abstract} This paper studies the shallow Ritz method for solving one-dimensional diffusion-reaction problems. The method is capable of improving the order of approximation for non-smooth problems. By following a similar approach to the one presented in \cite{cai2024fast}, we present a damped block Newton (dBN) method to achieve nearly optimal order of approximation. The dBN method optimizes the Ritz functional by alternating between the linear and non-linear parameters of the shallow ReLU neural network (NN). For diffusion-reaction problems, new difficulties arise: (1) for the linear parameters, the mass matrix is dense and even more ill-conditioned than the stiffness matrix, and (2) for the non-linear parameters, the Hessian matrix is dense and may be singular. This paper addresses these challenges, resulting in a dBN method with computational cost of ${\cal O}(n)$.

The ideas presented for diffusion-reaction problems can also be applied to least-squares approximation problems. For both applications, starting with the non-linear parameters as a uniform partition, numerical experiments show that the dBN method moves the mesh points to nearly optimal locations.

\end{abstract}



\section{Introduction}

Using neural networks (NNs) to solve partial differential equations (PDEs) has recently gained traction in scientific computing (see, e.g.,  \cite{BERG201828,  CAI2020109707, Do19, E_Yu_2018, RAISSI2019686, sirignano2018dgm}). In one dimension, the shallow ReLU NN generates a class of approximating functions equivalent to free knot splines (FKS). FKS can significantly improve the order of approximation for non-smooth functions and hence reduce the number of degrees of freedom dramatically (see \cite{Burchard74} and discussion in \cite{cai2024fast}). However, determining the optimal knot locations (the non-linear parameters of a shallow ReLU NN) leads to a complicated, computationally intensive non-convex optimization problem. Moreover, the ReLU activation function induces dense and ill-conditioned algebraic systems. 

For the one-dimensional diffusion problem, to realize optimal or nearly optimal order of the shallow Ritz approximation, we developed in \cite{cai2024fast} a damped block Newton (dBN) method that can efficiently move the uniformly distributed mesh points to nearly optimal locations. Due to the physical meaning of  the parameters of the output and hidden layers, the dBN method employs the commonly used outer-inner iterative method by alternating between updates of the linear and non-linear parameters (see, e.g., \cite{Golub1973, separable3, separable2} since the 1970s and \cite{Ainsworth2020PlateauPI, Ainsworth2022, pmlr-v107-cyr20a, park2022neuron} recently in the context of PDEs). This step is natural but resolves none of the essential difficulties mentioned above. The key contributions of the dBN method consist of (1) the exact inversion of the dense and ill-conditioned linear system and (2) one step of a damped Newton method applied to a {\bf reduced} non-linear system with $\mathcal{O}(n)$ computational cost per each iteration. Additionally, derivation of a modified Newton method for the non-linear parameters is highly nontrivial due to the facts that the ReLU activation function is not differentiable everywhere and that the corresponding Hessian matrix could be singular.

The purpose of this paper is to extend the dBN method to a broader class of problems in one dimension, while maintaining the efficiency previously achieved for diffusion problems (\cite{cai2024fast}). For diffusion-reaction problems as well as least-squares approximation, the mass matrix $M(\bb)$ arising from the NN approximation must be inverted per iteration for solving the linear parameter. Just as the stiffness matrix $A(\bb)$ in \cite{cai2024fast}, $M(\bb)$ depends on the non-linear parameter $\bb$. However, $M(\bb)$ is much more ill-conditioned than $A(\bb)$. Specifically, the condition numbers $\kappa(A(\bb))$ and $\kappa(M(\bb))$ are bounded, respectively, by $\mathcal{O}(nh_{\text{min}}^{-1})$ and $\mathcal{O}(nh_{\text{min}}^{-3})$ (see Lemma~4.3 in \cite{cai2024fast} and \cref{l:Cond}), where $n$ is the number of neurons and $h_{\text{min}}$ is the smallest distance between two neighboring breakpoints. For the ReLU NN approximation to non-smooth functions, $h_{\text{min}}$ could be much smaller than the uniform mesh size $1/n$ by several orders. To invert the dense matrices $M(\bb)$ and $A(\bb)+M(\bb)$ efficiently, we introduce their factorizations whose inversions can be computed in $\mathcal{O}(n)$ operations (see \cref{massMatrix_inv} and \cref{sum_inv}).

Along with the possible non-differentiability of the diffusion coefficient, the non-linear parameters for this broader class of problems present additional challenges. First, the optimality conditions for the non-linear parameters are no longer nearly decoupled non-linear algebraic systems (see \cref{opt_b_pde}), which imply dense Hessian matrices. Second, the Hessians could be singular and hence the Newton method is not applicable. To overcome the second difficulty, we remove some neurons, that are either unneeded or can be fixed, to obtain a reduced non-linear system. A neuron is not needed if its contribution to the current approximation is small, and a neuron can be fixed if it is at a nearly optimal location. Sufficient conditions for the derivation of the reduced system are introduced and guarantee the invertibility of the Hessian for the reduced non-linear system. To invert the dense Hessian matrix of the reduced system in $\mathcal{O}(n)$ operations, we utilize its special structure and the explicit formula for the inverse of the stiffness matrix in \cite{cai2024fast}. Hence, the overall computational cost per iteration of dBN remains ${\cal O}(n)$, consistent with the complexity established in \cite{cai2024fast}. 

To enhance the efficiency of moving breakpoints even further, we combine the dBN method with the adaptive neuron enhancement (ANE) method \cite{Liu_Cai_Chen_2022}. Numerical examples demonstrate the ability of the aforementioned methods to move the breakpoints quickly and efficiently and to outperform L-BFGS for select examples.

The paper is structured as follows. In \Cref{sec1}, we first describe the shallow Ritz discretization for diffusion-reaction problems. \Cref{sec2} presents the optimality conditions of the shallow Ritz discretization. An upper bound for the condition number of the mass matrix, along with a fast inversion strategy, is derived in \Cref{sec3}. The dBN method is introduced in \Cref{sec4}, and the adaptive version of dBN is described in \Cref{adaptivitySection}. The ideas presented in this paper make it possible to develop a dBN method for least-squares approximation problems, as explained in \Cref{ls_section}. Finally, numerical results are presented in \Cref{Section numerical exp}.

\section{Shallow Ritz Method for Diffusion-Reaction Problems} \label{sec1}Consider the following diffusion-reaction equation in one dimension:
\begin{equation}\label{pde}
    \left\{\begin{array}{lr}
        -(a(x)u^{\prime}(x))^{\prime} +r(x)u(x)=f(x), & \mbox{in }\, I=(0,1),\\[2mm]
        u(0)= \alpha,  \quad u(1)=\beta, &
    \end{array}\right.
\end{equation}
where the diffusion coefficient
$a(x)$, the reaction coefficient $r(x)$, and $f(x)$ are given real-valued functions defined on $I$. Assume that $a(x)\in L^\infty(I)$ and $r(x)\in L^\infty(I)$ are bounded below by the respective positive constant $a_0 > 0$ and non-negative constant $r_0\ge 0$ almost everywhere on $I$.

As in \cite{cai2024fast}, the modified Ritz formulation of problem \cref{pde} is to find $u \in H^{1}(I)\bigcap \{u(0)= \alpha\}$ such that
\begin{equation}\label{energy_functional}
J(u) = \min_{v \in H^{1}(I)\cap \{v(0)= \alpha\}}J(v),
\end{equation}
where the modified energy functional is given by 
\begin{equation}\label{modified_functional}
    J(v) = \frac{1}{2}\int_{0}^1a(x)(v^{\prime}(x))^2dx +\frac{1}{2}\int_{0}^1r(x)(v(x))^2dx -   \int_{0}^1f(x)v(x)dx + \frac{\gamma}{2}(v(1) - \beta)^2.
\end{equation}
Here, $\gamma > 0$ is a penalization constant.

The set of approximating functions generated by the shallow ReLU neural network with $n$ interior breakpoints over the domain $I$ is denoted by 
\begin{eqnarray*}\label{ShallowNN}
    {\cal M}_n(I)&=&\left\{c_{-1}+ c_0x+  \sum_{i=1}^{n}c_i\sigma(x-b_i) \, :\, c_i\in \R,\, 0 = b_0 <  b_1<\cdots<b_n<b_{n+1}=1\right\},  
\end{eqnarray*}
where  $\sigma(t) =\max\{0, t\}$ is the ReLU activation function. Then, the Ritz neural network approximation is to find $u_n\in {\cal M}_n(I) \cap \{u_n(0) = \alpha \}$ such that 
\begin{equation}\label{min_pde}
    J(u_n)=\min_{v\in\cM_n(I) \cap \{v(0)= \alpha\}}J(v).
\end{equation}

Next, we provide an error bound for the solution $u_n$ to \cref{min_pde}. The corresponding bilinear form of the modified energy functional is given by
\begin{equation*}
    a(u, v) := \int_0^1a(x)u^{\prime}(x)v^{\prime}(x)dx + \int_{0}^1r(x)u(x)v(x)dx +\gamma u(1)v(1)
\end{equation*}
for any $u,\,v\in H^1(I)$. 
Denote by $\|\cdot \|_{a}$ the induced norm of the bilinear form.

\begin{proposition}\label{l:error_estimate1}
Let $u$ and $u_n$ be solutions of problems \cref{energy_functional} and \cref{min_pde}, respectively. Then 
    \begin{equation}\label{LeaL}
        \|u-u_n\|_{a} \leq \sqrt{3}\inf_{v\in\cM_n(I) \cap \{v(0) = \alpha\}}\|u-v\|_{a} + 2\sqrt{2}\, \big|a(1)u^{\prime}(1)\big|\, \gamma^{-1/2}.
    \end{equation}
Moreover, if $\cM_n(I)$ has the following approximation property
\begin{equation}\label{app}
\inf_{v\in\cM_n(I)}\|u - v\|_{H^1(I)} \leq C(u)\, n^{-1},
\end{equation}
then there exists a constant $C$ depending on $u$ such that 
    \begin{equation}\label{Err}
        \|u - u_n\|_a \leq C \left(n^{-1}+\gamma^{-1/2}\right).
    \end{equation}
\end{proposition}
\begin{proof}
Inequality \cref{LeaL} may be proved in a similar fashion as that of Lemma 3.1 in \cite{cai2024fast}, and \cref{Err} is a direct consequence of \cref{LeaL} and \cref{app}.
\end{proof}


\section{Optimality Conditions}\label{sec2}This section derives systems of algebraic equations arising from the optimality conditions of \cref{min_pde}.

To this end, let 
\[
u_n(x) = \alpha+ c_0x+  \sum_{i=1}^{n}c_i\sigma(x-b_i)
\]
be a solution of \cref{min_pde} in ${\cal M}_{n}(I)\cap \{u_n(0) = \alpha \}$. Denote by $\bc = (c_0, \dots, c_n)^T$ and $\bb =(b_1, \dots, b_n)^T$ the respective linear and non-linear parameters. Then, the first-order optimality conditions yield
\begin{equation}\label{cps_ode}
    \nabla_{\bc} J\left(u_n\right)={\bf 0}
    \quad\mbox{and}\quad
    \nabla_{\bb} J\left(u_n\right)={\bf 0},
\end{equation}
 where $\nabla_{\bc}$ and $\nabla_{\bb}$ denote the gradients with respect to $\bc$ and $\bb$, respectively. 

 Let
 \begin{equation*}
   H(t) = \left\{\begin{array}{ll}
    1, & t >0,\\[2mm]
    \frac{1}{2}, & t = 0,\\[2mm]
    0, &t<0,
    \end{array}\right.\quad \mbox{and}\quad\delta(t) = \left\{\begin{array}{ll}
    +\infty, & t = 0,\\[2mm]
    0, & t \neq 0,   \end{array}\right. 
\end{equation*}
be the Heaviside (unit) step function and the Dirac delta function, respectively. Note that $H(t)=\sigma'(t)$ everywhere except at $t=0$, and similarly, $\delta(t)=\sigma''(t)$. Consider the mass matrix associated with the weight function $\omega$ given by 
\begin{equation}\label{massMatrix}
    M_\omega(\bb) = \big(m_{ij}\big) \quad\mbox{with }\, m_{ij}=\int_{0}^1\omega(x)\sigma(x-b_{i-1})\sigma(x -b_{j-1})dx 
\end{equation}
for $i, j = 1, \dots, n+1$, and the stiffness matrix associated with $\omega$ given by
\begin{equation*}\label{coeffMatrix}
    A_\omega(\bb) = \big(a_{ij}\big)\quad\mbox{with }\, a_{ij}=\int_{0}^1\omega(x)H(x-b_{i-k})H(x -b_{j-k})dx
\end{equation*} 
for $i, j = 1, \ldots, n+k$, where either $k = 0$ or $k=1$ depending on the context.

Denote the right-hand side vector by
\begin{equation*}
      \bff(\bb)=\big(f_{i}\big) \quad\mbox{with }\, f_{i}=  \int_0^1 \left(f(x)-\alpha\right) \sigma(x-b_{i-1})dx
\end{equation*}
for $i = 1, \ldots, n+1$ and let $\bd=\nabla_{\bc}u_n(1)$. By a similar derivation as in \cite{cai2024fast}, the first equation in \cref{cps_ode} becomes
    \begin{equation}\label{linear_eq}
       {\cal A}(\bb)\,\bc= {\cal F}(\bb) 
\end{equation}
where 
\begin{equation*}
    {\cal A}(\bb) = A_a(\bb)  + M_r(\bb) + \gamma \bd \bd^T \quad \mbox{and} \quad {\cal F}(\bb) = \bff(\bb) + \gamma (\beta - \alpha)\bd.
\end{equation*}
Comparing to (4.2) in \cite{cai2024fast}, the additional term $M_r(\bb)\bc$ in \cref{linear_eq} results from the reaction term.

Next, to derive the second equation in \cref{cps_ode}, for $j=1,\ldots,n$, let
\begin{equation}\label{q}
q_j = \int_{b_j}^1\bigl(f(x) -r(x)u_n(x)\bigr)dx  - a(b_j)u_n'(b_j), \quad \mbox{where} \quad
    u_n'(b_j): = \sum_{i=0}^{j-1}c_i+ \frac{c_j}{2}.
\end{equation}

In a similar fashion as in \cite{cai2024fast}, the second optimality condition in \cref{cps_ode} becomes
\begin{equation}\label{opt_b_pde}
    {\bf 0}=\nabla_{\bb} J\left(u_n\right)=\bD(\hat{\bc}) \left(\bq + \gamma(u_n(1) - \beta){\bf 1}\right),
\end{equation}
where $\bD(\hat{\bc})$ is a $n\times n$ diagonal matrix with the $i^{th}$ diagonal element $c_i$ and \[
\hat{\bc}= (c_1,\ldots, c_n)^T, \quad  \bq = (q_1, \ldots, q_n)^T, \quad\mbox{and}\quad {\bf 1} = (1, \ldots, 1)^T.
\]

\begin{remark}\label{c0_rmk}
     If $c_i = 0$, then there is no $i$-th equation in the optimality condition \cref{opt_b_pde}. Furthermore, the $i$-th neuron has no contribution to the approximation $u_n$. Such neurons can be removed or redistributed. 
\end{remark}

\section{Mass Matrix}\label{sec3}
This section studies the mass matrix resulting from a shallow ReLU neural network and the computation of its inversion.  

While the stiffness matrix $A_{a}(\bb)$ is dense, its inversion is a tri-diagonal matrix with an explicit algebraic formula (see \cite{cai2024fast}). This property holds for matrices with the following structure
\begin{equation}\label{specialMatrix}
    {\cal M}= 
    \begin{pmatrix}
         \alpha_1 \beta_1 & \alpha_1 \beta_2 & \alpha_1 \beta_3 &  \hdots & \alpha_1\beta_k\\
     \alpha_1 \beta_2 & \alpha_2\beta_2 & \alpha_2\beta_3 &   \hdots & \alpha_2 \beta_k\\
     \alpha_1\beta_3 & \alpha_2\beta_3 & \alpha_3\beta_3   & \hdots & \alpha_3\beta_k\\
     \vdots & \vdots & \vdots  &  \ddots & \vdots\\
     \alpha_1 \beta_k & \alpha_2\beta_k & \alpha_3\beta_k & \hdots & \alpha_k\beta_k
    \end{pmatrix},
\end{equation}
where $\alpha_i$ and $\beta_i$ are real numbers.

\begin{lemma}\label{inverseEspecialMatrix}
Assume that $\alpha_1 \neq 0$, $\beta_k \neq 0$, and $\alpha_{i+1}\beta_i - \alpha_{i}\beta_{i+1} \neq 0$
for $i = 1, \ldots, k-1$, then the matrix $\cM$ defined in \cref{specialMatrix}
is invertible. Moreover, its inverse $\cM^{-1}$ is symmetric and tri-diagonal with non-zero entries given by
\[
{\cal M}^{-1}_{ii}=\displaystyle\frac{\alpha_{i+1}\beta_{i-1} -\alpha_{i-1}\beta_{i+1}}{(\alpha_i\beta_{i-1} - \alpha_{i-1}\beta_i)(\alpha_{i+1}\beta_i - \alpha_i\beta_{i+1})}
\quad\mbox{and}\quad 
{\cal M}^{-1}_{i,i+1} = \displaystyle \frac{-1}{\alpha_{i+1}\beta_i -\alpha_i\beta_{i+1}} = {\cal M}^{-1}_{i+1, i},
\]
    where $\alpha_0 = \beta_{k+1} = 0$ and  $\alpha_{k+1} = \beta_{0} = 1$.
\end{lemma}

\begin{proof}
The assumptions imply that ${\cal M}^{-1}_{i,i}$ and ${\cal M}^{-1}_{i,i+1}={\cal M}^{-1}_{i+1,i}$ are well-defined. It is easy to verify that ${\cal M}{\cal M}^{-1} = I$.
\end{proof}

The stiffness matrix $A_a(\bb)$ has the same structure as ${\cal M}$ with 
\[
k= n, \quad \alpha_i= 1, \quad\mbox{and}\quad \beta_i=\int_{b_i}^1a(x)dx.
\] 
The mass matrix $M_r(\bb)$ is also dense due to the global support of the neurons, and its condition number is very large (see \cref{subsec2.0}). However, it can be factorized into matrices with structure as in \cref{specialMatrix}. This section concludes with the inverse formulas of the mass matrix, and its derivations are presented in \Cref{appendixA} algebraically and geometrically. By \cref{massMatrix_inv}, application of $M_r(\bb)^{-1}$ to any vector costs $\mathcal{O}(n)$ operations.

\subsection{Condition Number}\label{subsec2.0}

Let $h_i=b_{i+1}-b_{i}$ for $i=0,\ldots,n$, and set 
\[
h_{\text{max}}=\max\limits_{0\leq i\leq n} h_i \quad\mbox{and}\quad h_{\text{min}}=\min\limits_{0\leq i\leq n} h_i.
\]
It was shown in \cite{cai2024fast} that the condition number of $A_a(\bb)$ is bounded by $\mathcal{O}(nh_{\text{min}}^{-1})$ for $a(x) = 1$. For the mass matrix, it was established in \cite{hong2022activation} that the condition number is both lower and upper bounded by $\mathcal{O}(n^4)$ when the breakpoints are uniformly distributed over the domain $\Omega = [0,1]$. Building on this, the authors in \cite{caiMass} extended the analysis to arbitrary breakpoint distributions, relating the condition number to the minimum pairwise distance between breakpoints.

\begin{lemma}\label{l:Cond}
Let $r(x) = 1$, then the condition number of the mass matrix $M_r(\bb)$ is bounded by $\mathcal{O}\left(n/h_{\textup{min}}^3\right)$.
\end{lemma}

\begin{proof}
For any vector $\bxi =(\xi_0, \xi_1,\ldots,\xi_n)^T\in \mathbb{R}^{n+1}$, denote its magnitude by $\big|\bxi\big|=\left(\sum\limits_{i=0}^n \xi^2_i\right)^{1/2}$. By the Cauchy-Schwarz inequality and the fact that $\sigma(x-b_j)=0$ for $x\leq b_j$,
we have
    \begin{align}\label{cond_upper}
    \begin{split}
    \bxi^T M_r(\bb) \bxi 
    &= \int_0^1 \left( \sum_{i=0}^{n}\xi_i\sigma(x-b_i) \right)^2 \, dx 
    \leq |\bxi|^2 \int_{0}^1 \left(\sum_{i=0}^n \sigma(x-b_i)^2\right) \, dx  \\
    &= |\bxi|^2 \sum_{j=0}^n \int_{b_j}^{b_{j+1}}\sum_{i=1}^{j}\sigma(x-b_i)^2 \,dx
    = \frac{|\bxi|^2}{3} \sum_{j=0}^n \sum_{i=0}^j \left\{(b_{j+1}-b_i)^3-(b_j-b_i)^3\right\}\\
    &= \frac{|\bxi|^2}{3}\sum_{i=0}^n (b_{n+1}-b_i)^3
    = \frac{|\bxi|^2}{3}\sum_{i=0}^n (1-b_i)^3\leq \frac{(n+1)}{3}\, |\bxi|^2.
    \end{split}
    \end{align}
    
To estimate the lower bound of $\bxi^T M_r(\bb) \bxi$, let 
\[
\tau_i(x)=\sum\limits_{j=0}^i \xi_j \sigma(x-b_j) \quad\mbox{and}\quad a_{i-1}= \tau_i(b_i)=\sum\limits_{j=0}^{i}\xi_j(b_i-b_j),
\]
for $i = 0, \dots, n+1$, where $a_{-1} = 0$. Then $\tau_i\left(\frac{b_{i+1}+b_i}{2}\right)=\dfrac{a_{i-1}+a_{i}}{2}$. Since 
$\tau_i^2(x)$ is a quadratic function in each sub-interval $[b_i,b_{i+1}]$, Simpson's Rule implies
    \begin{align*}\label{cond_lower}
    \begin{split}
       \bxi^T M_r(\bb) \bxi 
    &=\sum_{i=0}^n\int_{b_i}^{b_{i+1}}\tau^2_i(x)\,dx
    =\frac16\sum_{i=0}^n h_{i}\left[ \tau_i^2(b_i) + 4\tau_i^2\left(\frac{b_{i+1}+b_i}{2}\right) + \tau_i^2(b_{i+1}) \right]\\
    &= \frac16 \sum_{i=0}^n h_{i}\left[a_{i-1}^2+(a_{i-1}+a_{i})^2+a_{i}^2\right]
    \geq \frac{1}{6}h_{\text{min}}\lvert{\bf a}\rvert^2,
    \end{split}
    \end{align*}
where $\ba=(a_0, a_1,\ldots,a_n)^T$. It is easy to see that 
$\bxi= Q \ba$, where $Q$ is a $(n+1)$-order lower tri-diagonal matrix given by
\begin{equation*}\label{Q}
    Q= \begin{pmatrix}
      ~\frac{1}{h_0} & 0 &  0  & \dots & 0 & 0\\ 
    -\left(\frac{1}{h_0}+\frac{1}{h_1}\right) &   \frac{1}{h_1} &  0  & \dots & 0 &  0  \\[2mm]
      \frac{1}{h_1}   & -\left(\frac{1}{h_1}+\frac{1}{h_2}\right) & \frac{1}{h_2} & \dots & 0 & 0  \\
      \vdots & \vdots & \vdots & \ddots & \vdots & \vdots\\ 
0 & 0 & 0 & \hdots & \frac{1}{h_{n-1}} & 0\\ 
      0 & 0 & 0 & \hdots & -\left(\frac{1}{h_{n-1}}+\frac{1}{h_{n}}\right) & \frac{1}{h_{n}}\\
  \end{pmatrix}.
\end{equation*}
It is easy to verify that $Q$ has spectral norm bounded by
\begin{equation*}
    \lVert Q\rVert_2\leq \sqrt{\lVert Q\rVert_1\lVert Q\rVert_{\infty}}\leq 4h_{\text{min}}^{-1}.
\end{equation*}
Hence,
\[
 \bxi^T M_\omega(\bb) \bxi \ge  
 \frac16 h_{\text{min}} \frac{\lvert{\bxi} \rvert^2}{\lVert Q\rVert_2^2} \ge  \frac{1}{96} h^3_{\text{min}} \lvert{\bxi}\rvert^2,
\]
which, together with the upper bound in \cref{cond_upper}, implies the validity of the lemma. 
\end{proof}

\begin{lemma}\label{Mr}
Under the assumption on the reaction coefficient $r(x)$,
the condition number of the mass matrix $M_r(\bb)$ is bounded by $\mathcal{O}\left(n\,r^{-1}_0 h_{\textup{min}}^{-3}\right)$. 
\end{lemma}

\begin{proof}
Since $r\in L^\infty(I)$ and $r(x)\geq r_0$ almost everywhere, in a similar fashion as the proof of \cref{l:Cond}, we have 
\[
\dfrac16 r_0 h_{\text{min}}^3 \lvert{\bxi}\rvert^2 \leq \bxi^T M_r(\bb) \bxi \leq C\, (n+1)  \lvert{\bxi}\rvert^2,
\]
which implies the validity of the lemma.
\end{proof}

Whereas the mass matrix associated with the ReLU neural network is very ill-conditioned, it is well-known that the mass matrix for the finite element (FE) method is much better conditioned (see \cite{Fried_1973} for example). The following \cref{l: FE_cond} reiterates the result in \cite{Fried_1973} with an alternate proof in a similar fashion as that of \cref{l:Cond}. To this end, for $i=1,\ldots, n$, denote the standard hat basis functions by
  \begin{equation*}\label{hat_fxns}
        \varphi_i(x)=\left\{\begin{array}{ll}
    (x-b_{i-1})/h_{i},   & x\in (b_{i-1},b_i), \\[2mm]
    (b_{i+1}-x)/h_{i+1},  & x\in (b_{i},b_{i+1}), \\[2mm]
    0, & \text{otherwise}.
    \end{array}\right.
    \end{equation*}
Let $\bphi=\left(\varphi_1,\ldots,\varphi_{n}\right)^T$, then the corresponding mass matrix is given by 
    \begin{equation*}\label{fe_mass}
        \Tilde{M}(\bb) =\int_0^1 \bphi\bphi^Tdx.
    \end{equation*}


\begin{lemma}\label{l: FE_cond}
The condition number of the finite element mass matrix $\Tilde{M}(\bb)$ is bounded by\\ $\mathcal{O}(h_{\textup{max}}/h_{\textup{min}})$.
\end{lemma}

\begin{proof}
For any vector $\bxi =(\xi_1,\ldots,\xi_{n})^T\in \mathbb{R}^{n}$, 
in a similar fashion as that of \cref{l:Cond}, 
we get the equality
    \begin{align*}\label{FE_lower_bd}
    \begin{split}
        \bxi^T \Tilde{M}(\bb) \bxi &= \sum_{j=0}^{n} \int_{b_{j}}^{b_{j+1}}\left( \xi_j\varphi_j+\xi_{j+1}\varphi_{j+1}\right)^2\,dx 
        =\sum_{j=0}^{n} \frac{h_{j}}{6}\left[ \xi_j^2+(\xi_j+\xi_{j+1})^2+\xi_{j+1}^2 \right]\\
    \end{split}
    \end{align*}
    with $\varphi_{0}(x) = \varphi_{n+1}(x) = \xi_0= \xi_{n+1}=0$, which leads to the inequalities
    \begin{equation*}
        \frac{1}{6}h_{\text{min}}|\bxi|^2 \leq \bxi^T \Tilde{M}(\bb) \bxi\leq \frac{2}{3}h_{\text{max}}|\bxi|^2.
    \end{equation*}
This completes the proof of the lemma.
\end{proof}
\begin{rem}
It is well-known {\em \cite{Loan_Golub}} that classic iterative solvers such as the Richardson (i.e., the method of gradient descent), Jacobi, Gauss-Seidel, etc. for a system of linear equations converge slowly if the condition number $\kappa(A)$ of the corresponding matrix is very large. The convergence rate of those methods is at most $\frac{\kappa(A)-1}{\kappa(A)+1}$.
\end{rem} 
\subsection{Inversion of the Mass Matrix}\label{mass_inversion}
Despite the ill-conditioning of the mass matrix $M_r(\bb)$, it can still be efficiently inverted. In fact, the following lemma states that $M_r(\bb)$ can be factorized into a product of matrices that are easy to invert. The factorization can be derived by either an algebraic (see \cref{appendixA_algebrajc}) or a geometric (see \cref{appendixA_geometric}) approach.

\begin{lemma}
The mass matrix $M_r(\bb)$ has the following factorization
    \begin{equation}\label{massMatrix_inv}
        M_r(\bb) = T_1^{-1}T_2T_3^{-1},
    \end{equation}
where $T_{1}$, $T_{2}$, and $T_3$ are tri-diagonal matrices given in \cref{inv-M} and \cref{fac}.
\end{lemma}
 

 \begin{remark}
Additionally, for the given $T_1, T_3$ in \cref{fac}, the sum $A_a(\bb) + M_r(\bb)$ has the following factorization \begin{equation}\label{sum_inv}
      A_a(\bb)+ M_r(\bb)  = T_1^{-1}T_4 T_3^{-1},
 \end{equation}
 where $T_4$ is a tri-diagonal matrix given in \textit{\textit{\cref{rmk: coeff_factorization}}}.
 \end{remark}

 \begin{remark}
     The Sherman-Morrison formula and \cref{sum_inv} imply that the system of linear equations in \cref{linear_eq} can be solved in ${\cal O}(n)$ operations.
 \end{remark}

 \section{A Damped Block Newton Method}\label{sec4} As mentioned in \cite{cai2024fast}, the system of non-linear algebraic  equations in \cref{opt_b_pde} can be solved using the Variable Projection (VarPro) method \cite{Golub1973}. In the VarPro approach, the solution $\bc$ of \cref{linear_eq} is substituted into \cref{opt_b_pde}, resulting in a reduced non-linear system for $\bb$. However, this reduced system may exhibit a more complex structure.

In this section, we describe a damped block Newton method for solving the minimization problem in \cref{min_pde}. This method employs the block Gauss-Seidel method as an outer iteration for the linear and non-linear parameters. Per each outer iteration, the linear and the non-linear parameters are updated by exact inversion and one step of a damped Newton method, respectively.

To solve \cref{linear_eq} for the linear parameter, the direct inversion is addressed in \cref{sec3}. Now to apply Newton's method to the optimality condition \cref{opt_b_pde}, we first must derive the Hessian for the non-linear parameter and discuss its invertibility. To that end, assume that $a(x)$ is differentiable at $b_i$ for  each $i = 1, \dots, n$, and define
\begin{equation*}
    g_i = \frac{\partial}{\partial b_i}q_i = r(b_i)u_n(b_i) - f(b_i)-  a^\prime(b_i)u_n'(b_i),
\end{equation*} 
where $q_i$ is defined in \cref{q}. Denote by $\bD(\bg) = \text{diag} (g_1, \dots, g_n)$
the diagonal matrix with the $i^{th}$ diagonal element $g(b_i)$. A formula for the Hessian matrix $\nabla_{\bb}^2J$ is given in the next lemma.

\begin{lemma}\label{hessian}
Assume that $c_i \neq 0$ and that $a(x)$ is differentiable at $x=b_i$ for all $i =1, \dots, n$. Then the Hessian matrix $\nabla^2_{\bb}J(u_n)$ has the form
\begin{equation}\label{H_DR}
    \mathbf{H}(\bc, \bb) = \bD(\hat{\bc})\bigl(\bD(\bg)
    + A_{r}(\bb)\bD(\hat{\bc}) + \gamma {\bf 1}\hat{\bc}^T\bigr).
    \end{equation}
    \end{lemma}
\begin{proof}
Equation \cref{H_DR} can be derived in a similar fashion as that of Lemma 5.2 in \cite{cai2024fast}. The only difference here is the additional reaction term in \cref{modified_functional}. For that term, the computations shown in Lemma~4.1 from \cite{Cai24GN} can be used to obtain the second-order derivatives with respect to $\bb$.
\end{proof}

The following lemma gives a sufficient condition for the invertibility of the Hessian $\bH(\bc, \bb)$.

\begin{lemma}\label{sufff}
Let $\tilde{h}_i=min\{h_{i-1},h_{i}\}$. If $c_i \neq 0$, $a(x)$ is differentiable at $b_i$ and
    \begin{equation}\label{inv_cond}
        \frac{g_i}{c_i}+\frac{r_0\tilde{h}_i}{4}>0 \quad\mbox{for all}\quad i=1,\dots,n,
    \end{equation}
    then $\bD(\bg)\bD(\hat{\bc})^{-1}+A_r(\bb)$ is positive definite.
\end{lemma}

\begin{proof}
For any vector $\bxi =(\xi_1,\ldots,\xi_n)^T\in \mathbb{R}^{n}$, a similar argument to the proof of Lemma 4.3 in \cite{cai2024fast} shows that
    \[
        \bxi^T A_r(\bb)\bxi\geq  \frac{r_0}{4}\sum_{j=1}^n \tilde{h}_j\xi_j^2, 
    \]
which implies that $ \bxi^T \left(\bD(\bg)\bD(\hat{\bc})^{-1}+ A_r(\bb)\right)\bxi\geq  \sum\limits_{j=1}^n \left(\frac{g_i}{c_i}+\frac{r_0}{4}\tilde{h}_j\right)\xi_j^2$. Now, the lemma is a direct consequence of the assumption in \cref{inv_cond}.
\end{proof}

\begin{corollary}\label{cor}
     The matrix ${\bf H}(\bc,\bb)$ is invertible under the conditions of \cref{sufff}.
\end{corollary}

\begin{remark}
 Condition \cref{inv_cond} is sufficient but not necessary. Therefore, for implementation purposes, we consider the relaxed condition
\begin{equation}\label{relaxed_conv} \frac{g_i}{c_i} + \tau_2^{-1}\frac{r_0\tilde{h}_i}{4}   > 0 \quad \text{for all} \quad i = 1, \dots, n, \end{equation}
for some parameter $0< \tau_2 \leq 1$. In the numerical examples presented in \cref{Section numerical exp}, setting $\tau_2 = 1$ was not necessary to ensure invertibility, so a value $\tau_2 < 1$ was used instead.
\end{remark}

\begin{remark}
    The action of ${\bf H}(\bc,\bb)^{-1}$ applied to any vector can be computed in ${\cal O}(n)$ operations.  This is due to the Sherman-Morrison formula and the facts that 
    \[  \bD(\bg)+A_r(\bb)\bD(\hat{\bc})=\left(\bD(\bg)\bD(\hat{\bc})^{-1}A_r(\bb)^{-1}+I\right)A_r(\bb)\bD(\hat{\bc}),
    \]
    and that $A_r(\bb)^{-1}$ is tri-diagonal.  
\end{remark} 

In cases where $a'(b_i)$ is not defined for some $i \in \{1, \dots, n\}$, the breakpoint $b_i$ lies on the interface and should be fixed without further update. If $c_i = 0$, then the breakpoint $b_i$ could either be removed or redistributed (see \cref{c0_rmk}). 

Overall, we can update the non-linear parameter with a Newton step in $\mathcal{O}(n)$ operations granted that the Hessian $\mathbf{H}(\bc, \bb)$ is invertible and the optimality condition \cref{opt_b_pde} is well-posed. Thus, for the implementation of the dBN method, we set aside the neurons which violate these conditions and construct a reduced system with said neurons removed. To that end, let $0 < \tau_1$ and $0<\tau_2 \leq 1$, and define the set
\[S_1 = \left\{i\in\{1,\dots,n\}: \left\lvert c_i\right\rvert<\tau_1 \text{ or } b_i\notin I\right\}\]
of indices for the non-contributing neurons, and the set
\[
S_2=\left\{i\in\{1,\dots,n\}\backslash S_1: \frac{g_i}{c_i}+ \tau_2^{-1}\frac{r_0\tilde{h}_i}{4} \leq  0 \quad\mbox{or}\quad a'(b_i)\quad\mbox{DNE}  \right\}
\]
of indices for which the corresponding neurons do not satisfy the invertibility condition \cref{relaxed_conv} or belong to the interface. Then
\begin{equation*}
    S = \{1, \dots, n\} \backslash (S_1 \cup S_2)
\end{equation*}
is the set of indices for the neurons which remain in the system.


Next, given a vector $\bv \in \mathbb{R}^n$, we denote by $\bv_{S} \in \mathbb{R}^{n-\left\lvert S^c\right\rvert}$ as the vector obtained by removing the $i-th$ entry from $\bv$ for every $i \in S^c=\{1,\dots, n\} \backslash S$. Similarly, for a matrix $\bB \in \mathbb{R}^{n\times n}$, define $\bB_{S} \in \mathbb{R}^{(n-\lvert S^c\rvert)\times (n-\lvert S^c\rvert)}$ as the matrix obtained by removing the $i-th$ row and column of $\bB$ for each $i \in S^c$.

Finally, define the \textit{reduced} search direction vector for the Newton step as
\begin{equation}\label{reduce_Newt}
    \bd_{R}(\bc, \bb)=  -\biggl[\bigl(\bD(\bg)
    + A_{r}(\bb)\bD(\hat{\bc}) + \gamma {\bf 1}\hat{\bc}^T\bigr)_{S}\biggr]^{-1}\bigg(\bq + \gamma(u_n(1) - \beta){\bf 1}\bigg)_{S}.
\end{equation}


We are now ready to describe the dBN method (see \cref{alg:dBN} for a pseudocode). Given  prescribed tolerances $ 0 <\tau_1$ and $0 <\tau_2 \leq 1$, let $\bb^{(k)}$ be the previous iterate, then the current iterate $\br^{(k+1)} =\begin{pmatrix}  \bc^{(k+1)} \\ \bb^{(k+1)} \end{pmatrix}$ is computed as follows:
\begin{itemize}
    \item[(i)] \textit{Compute the linear parameters} 
    \begin{equation*}
    \bc^{(k+1)} = {\cal A}\left(\bb^{(k)}\right)^{-1} {\cal F}\left(\bb^{(k)}\right).
    \end{equation*}
    \item[(ii)] \textit{Compute the search direction $\mathbf{p}^{(k)}=\left(p_1^{(k)},\dots,p_n^{(k)}\right)^T$} by
    \begin{equation}\label{directionVector}
    \bigl(p^{(k)}_i \bigr)_{i \in S} = d_R(\bc, \bb) \quad \mbox{and} \quad \bigl(p^{(k)}_i\bigr)_{i \notin S} = {\bf 0},
    \end{equation}
    where $d_R(\bc, \bb)$ is the \textit{reduced} Newton's direction vector defined in \cref{reduce_Newt}.
    \item [(iii)] \textit{Calculate the stepsize} $\eta_k$ by performing a one-dimensional optimization
    \begin{equation*}
        \eta_k = \arg\!\!\min\limits_{\eta \in \R^+}J\left(u_{n}\left(x; \bc^{(k+1)}, \bb^{(k)} + \eta\mathbf{p}^{(k)}\right)\right).
    \end{equation*}
    \item [(iv)] \textit{Compute the non-linear parameters}
    \begin{equation*}
        \bb^{(k+1)} = \bb^{(k)} +\eta_k\mathbf{p}^{(k)}.
    \end{equation*}
    \item [(v)] \textit{Redistribute non-contributing breakpoints} $b_i^{(k+1)}$ for all $i \in S_1$ and sort $\bb^{(k+1)}$.   
\end{itemize}


\begin{remark}\label{rmk_redist}
    The redistribution implemented for the numerical results shown in \cref{Section numerical exp} in step (v) was carried out as follows: For a neuron $b_l^{(k+1)}$ satisfying $l \in S_1$, we set 
    \[
    b_l^{(k+1)} \leftarrow\frac{b^{(k+1)}_{m-1}+b^{(k+1)}_{m}}{2},
    \]
    where $m\in\{1,\dots,n+1\}$ is an integer chosen uniformly at random.
\end{remark}

\begin{algorithm}
    \caption{A damped block Newton (dBN) method for  \cref{min_pde}} \label{alg:dBN} 
    \begin{algorithmic}
    \REQUIRE{Initial network parameters $\bb^{(0)}$}
    \ENSURE{Network parameters $\bc$, $\bb$}
    \FOR{$k=0,1\ldots$}
    \STATE $\triangleright$ \textit{Linear parameters}
    \STATE{$\bc^{(k+1)}\leftarrow {\cal A}\left(\bb^{(k)}\right)^{-1} {\cal F}\left(\bb^{(k)}\right)$} 
    \STATE $\triangleright$
    \textit{Non-linear parameters}
    \STATE{{\textit{Compute the search direction} }$\mathbf{p}^{(k)}$ \textit{as in}} \cref{directionVector}
    \STATE{$\eta_k \leftarrow \arg\!\!\min\limits_{\eta \in \R^+}J(u_n(x; \bc^{(k+1)}, \bb^{(k)} + \eta\mathbf{p}^{(k)}))$}
    \STATE{$\bb^{(k+1)} \leftarrow \bb^{(k)} + \eta_k\mathbf{p}^{(k)}$}
    \STATE $\triangleright$ \textit{Redistribute non-contributing neurons and sort $\bb^{(k+1)}$}
\ENDFOR
    \end{algorithmic}
\end{algorithm}

\begin{remark}    \label{initialize} The minimization problem in \cref{min_pde} is non-convex and contains multiple local minima. Therefore, obtaining an initial approximation sufficiently close to the desired minimum is crucial. Since the non-linear parameters $\bb$ correspond to the breakpoints that partition the interval $I=[0,1]$, we initialize them following the approach in \cite{Cai2021linear, cai2024fast, Liu_Cai_2022}. Specifically, the non-linear parameters $\bb$ are set as a uniform partition of $I$,  while the linear parameters $\bc$ are determined as the solution to the linear system in \cref{linear_eq} for this uniform mesh.
\end{remark}




\subsection{An Adaptivity Scheme}\label{adaptivitySection}
For a fixed number of neurons, the dBN method for the diffusion-reaction equation moves the initial uniformly distributed breakpoints very efficiently to nearly optimal locations as shown in \cref{Section numerical exp}. However, it was shown in \cite{cai2024fast} that introducing adaptivity results in a more optimal convergence rate.

In fact, the adaptive neuron enhancement (ANE) method \cite{Liu_Cai_2022, Liu_Cai_Chen_2022} was employed in \cite{cai2024fast}. 
The ANE method starts with a relatively small neural network and adaptively adds new neurons based on the previous approximation. At each adaptive step, the dBN method numerically solves the minimization problem in \cref{min_pde}. Subsequently, the newly added neurons are initialized based on regions where the previous approximation is not accurate.

In the ANE method, new neurons are added according to a marking strategy. Below, we explain the local indicators and the marking strategy used in this paper.

Let $\mathcal{K} = [c, d] \subseteq [0, 1]$ be a subinterval, then a modified local indicator of the ZZ type on $\mathcal{K}$ (see, e.g., \cite{CaZh:09}) is defined by
\[\xi_{\mathcal{K}}^2 =  \lVert a^{-1/2}\left(G(a u_n')-a u_n'\right) \rVert_{L^2(\mathcal{K})}^2 +  (d-c)^2\lVert -G'(a^2 u_n') + u_n -  f \rVert_{L^2(\mathcal{K})}^2,\]
where $G(v)$ is the projection of $v$ onto the space of the continuous piecewise linear functions. The corresponding relative error estimator is defined by $\xi = \xi_{I}/|u_n|_{H^1(I)}$.

Let $u_n \in {\cal M}_n(I)$ be a NN with the breakpoints 
\[
0=b_0 <b_1< \ldots < b_n<b_{n+1}=1.
\]
Define $\mathcal{K}^i = [b_{i},b_{i+1}]$, so that the collection ${\cal K}_n = \left\{ \mathcal{K}^i \right\}_{i = 0}^{n}$ defines a partition of the interval $[0, 1]$.
To refine this partition, we define $\widehat{\mathcal{K}}_n\subset \mathcal{K}_n$ by using the following average marking strategy:
\begin{equation}\label{marking}
    \widehat{\mathcal{K}}_n=\left\{K\in\mathcal{K}_n: \xi_{\mathcal{K}}\geq\frac{1}{\#\mathcal{K}_n}\sum_{{\cal K}\in\mathcal{K}_n}\xi_{\mathcal{K}}\right\},
\end{equation}
where $\#\mathcal{K}_n$ is the number of elements in $\mathcal{K}_n$.The AdBN  method is described in \cref{alg:adBN}.

\begin{algorithm}
    \caption{Adaptive damped block Newton (AdBN) method}\label{alg:adBN}
    \begin{algorithmic}
        \REQUIRE Initial number of neurons $n_0$, 
        parameters $a(x)$, $r(x)$ $f(x)$, $\alpha$, and $\beta$, tolerance $\epsilon$,
        \STATE (1)  Compute an approximation to the solution $u_n$ of the optimization problem in \cref{min_pde} by the dBN method; 
        \STATE (2) Compute the estimator $\xi_n=\left(\sum\limits_{K\in\mathcal{K}}\xi_K^2\right)^{1/2}/|u_n|_{H^1(I)}$;
        \STATE (3) If $\xi_n \leq \epsilon$, then stop; otherwise go to step (4);
        \STATE (4) Mark elements in $\widehat{\mathcal{K}}_n$ and denote by $\#\widehat{\mathcal{K}}_n$ the number of elements in $\widehat{\mathcal{K}}_n$;
        \STATE (5) Add $\#\widehat{\mathcal{K}}_n$ new neurons to the network and initialize them at the midpoints of elements in $\widehat{\mathcal{K}}_n$, then go to step (1)
    \end{algorithmic}
\end{algorithm}

\section{Least-Squares Approximation}\label{ls_section} 


Given a function $u(x)\in L^2(I)$, the best least-squares approximation to $u$ in $\cM_n(I)$ is to find $u_n \in \cM_n(I)\cap\{u_n(0)=u(0)\}$ such that
\begin{equation}\label{min}
    J(u_n)=\min_{v\in\cM_n(I)\bigcap \{v(0)=u(0)\}}J(v),
\end{equation}
where $J(v)$ is the weighted continuous least-squares loss functional given by
\begin{equation}\label{ls_fxnl}
    J(v) = \dfrac{1}{2}
\int_0^1 r(x)\left(v(x) - u(x)\right)^2dx=\frac{1}{2}\lvert\lvert v-u\rvert\rvert^2_{r}. 
\end{equation}


Let $u_n(x)=u(0) + \sum\limits_{i=0}^{n}c_i\sigma(x - b_i)\in\cM_n(I)$ be a solution of \cref{min}. Then, the optimality conditions become 
\begin{equation}\label{lineaEqDF}
    {\bf 0}=\nabla_{\bc} J\left(u_n\right)=M_{r}(\bb)\,\bc -  \bff(\bb) \quad\mbox{and}\quad {\bf 0}=\nabla_{\bb} J\left(u_n\right)= \bD(\hat{\bc})\br,
\end{equation}
where $M_{r}(\bb)$ is the mass matrix defined in \cref{massMatrix}, $\bD(\hat{\bc}) = \text{diag}(c_1, \dots, c_n)$, and $\bff(\bb)$ and $\br$ are given respectively by
\begin{align*}
      \bff(\bb)=\big(f_{i}\big)_{(n+1)\times 1} &\quad\mbox{with }\, f_{i}=  \int_0^1 r(x)\left(u(x)-u(0)\right)\sigma(x-b_{i-1})dx \\[0mm]
     \mbox{and}\quad \br =\big(r_{i}\big)_{n\times 1} &\quad\mbox{with }\,
      r_i=- \int_{b_i}^{1} r(x)(u_n(x) -u(x)) dx.
\end{align*}
Let $w_i = r(b_i)\bigl(u_n(b_i) - u(b_i)\bigr)$ for $i =1 ,\ldots, n$. In one dimension, Lemma~4.1 in \cite{Cai24GN} implies that the corresponding Hessian matrix is of the form
\begin{equation}\label{Hessianls}
\nabla_{\bb}^2 J(u_n)\equiv {\bf H}(\bc, \bb) = \bD(\hat{\bc})\bigl(\bD(\bw) + A_r(\bb)\bD(\hat{\bc})\bigr),
\end{equation}
where $\bD(\bw)=\text{diag}(w_1, \dots, w_n)$ is a diagonal matrix. 
Based on the optimality condition in \cref{lineaEqDF} and the Hessian matrix in \cref{Hessianls}, we can then design the dBN and the AdBN methods in a similar fashion as in \Cref{sec3} and \Cref{sec4}, respectively.

We conclude this section by discussing the numerical integration of the functionals in \cref{modified_functional} and \cref{ls_fxnl}. In one dimension, we may calculate the integrals analytically, as shown in several examples in \Cref{Section numerical exp}. In practice, integrals are often computed by various numerical integration methods. For an integrand (given as in \cref{ls_fxnl} or unknown as in \cref{modified_functional}), choosing an accurate quadrature rule with the least computational cost is highly non-trivial. An efficient and effective method for achieving this goal is the so-called adaptive quadrature (see, e.g., \cite{rice1975metalgorithm}), which was used in \cite{LiCaRa23} for the deep Ritz method in linear elasticity. 

An accurate but possibly inefficient numerical quadrature is, for example, the midpoint rule on a very fine uniform {\it integration} mesh. Specifically, let $\left\{x_i=i/m\right\}_{i=0}^m$ be a partition of the interval $I=[0,1]$ with sufficiently large $m$. Then the least-squares functional $J(v)$ is approximated by the following discrete least-squares functional
\begin{align*}
    J_{{m}}(v) = \dfrac{1}{2}
\sum_{i=1}^m r(x_{i-1/2})\left(v(x_{i-1/2}) - u(x_{i-1/2})\right)^2=\frac{1}{2}\lvert\lvert v-u\rvert\rvert^2_{r,m}. 
\end{align*}
The accuracy of the approximation depends on both the complexity of the function $u(x)$ and the size of $m$. Now, the discrete least-squares approximation to $u$ in $\cM_n(I)$ is to find $u_n \in \cM_n(I)\cap\{u_n(0)=u(0)\}$ such that
\begin{equation}\label{discrete_min}
    J_{{m}}(u_n)=\min_{v\in\cM_n(I)\bigcap \{v(0)=u(0)\}}J_{{m}}(v).
\end{equation}
The dBN developed for the continuous optimization problem in \cref{min} can be applied directly to the discrete optimization in \cref{discrete_min} by replacing integration with summation when forming the matrices $M_r(\bb)$ and $A_r(\bb)$. Since $\|\cdot\|_{r,m}$ defines a norm in $\mathbb{R}^m$, then it is easy to see that $M_r(\bb)$ and $A_{r}(\bb)$ are positive definite (see, e.g., Lemma~4.1 of \cite{Cai24GN}).

Treating $\left\{x_{i-1/2},u(x_{i-1/2})\right\}_{i=1}^m$ as a training set, the minimization problem in \cref{discrete_min} is similar to the standard machine learning problem in data science. Hence, the dBN method can be applied directly and viewed as an efficient training algorithm.

\section{Numerical Experiments} \label{Section numerical exp}

This section first presents numerical results using the dBN method to solve \cref{min}. Subsequently, results for the dBN and AdBN methods applied to \cref{pde} are shown in \cref{sectionex1} and \cref{sing_perturb}. The parameters $\tau_1$ and $\tau_2$ were set to $10^{-10}$ and $10^{-2}$, respectively. For diffusion-reaction problems, the penalization parameter $\gamma$ was set to $10^4$. In the AdBN method, a refinement occurred when the difference in relative residuals between two consecutive iterates was less than $10^{-6}$.

For each test problem of the diffusion-reaction equation, let $u$ and $u_n$ be the exact solution and its approximation in $\mathcal{M}_n(I)$, respectively. Denote the relative error by
\begin{equation}\label{seminorm_error}
    e_n = \frac{|u - u_n|_{H^1(I)}}{|u|_{H^1(I)}}.
\end{equation}

\subsection{Least-Squares Problem}\label{num ex 2}
\begin{figure}[htb!]
    \centering
    \subfigure[Loss versus number of iterations with 24 neurons. Final losses: L-BFGS - $6.78 \times 10^{-6}$, dBN - $4.16 \times 10^{-7}$]
    {\includegraphics[width=0.45\textwidth]{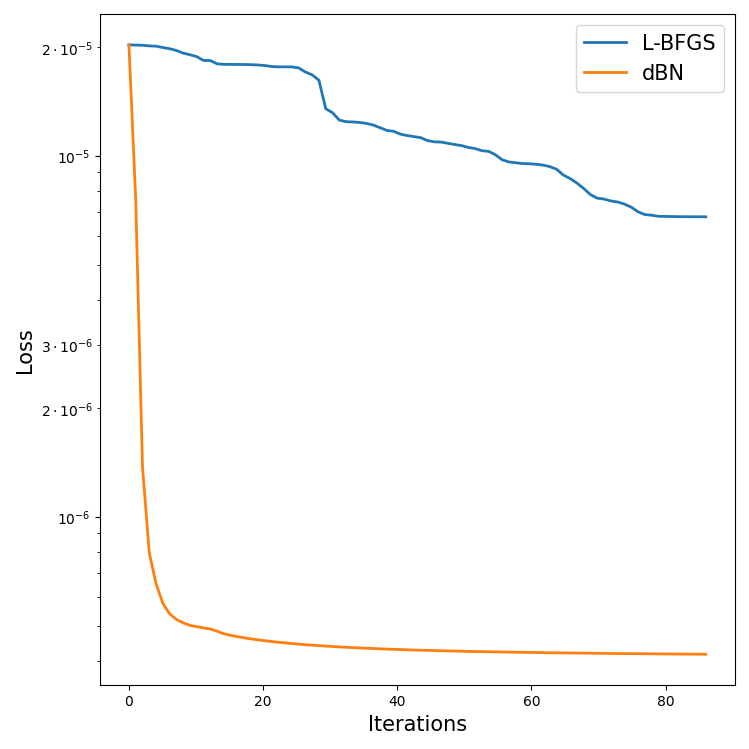}}
    \hspace{2em}
    \subfigure[Loss versus number of iterations with 48 neurons. Final losses: L-BFGS - $2.94 \times 10^{-6}$, dBN - $1.07 \times 10^{-7}$]{
        \includegraphics[width=0.45\textwidth]{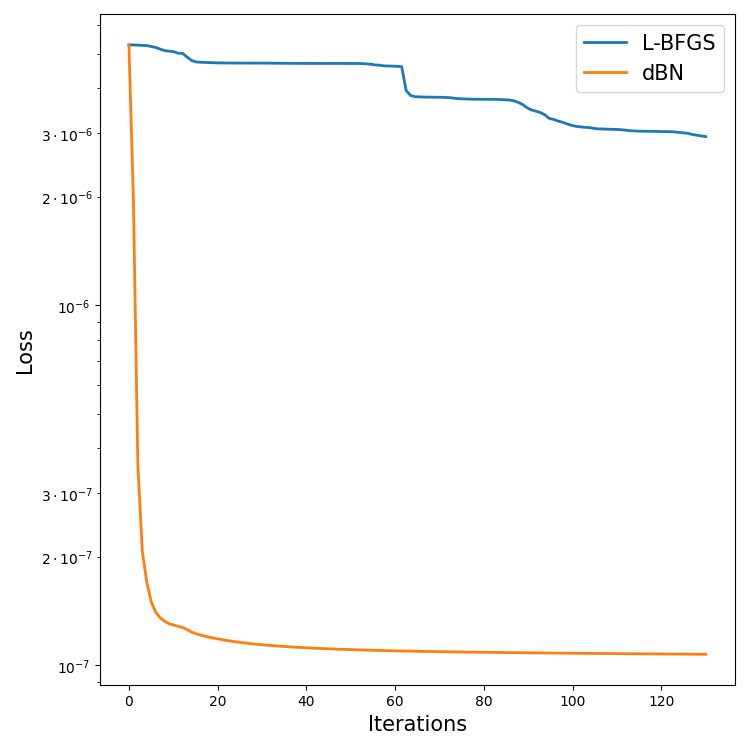}}
    \caption{Comparison between L-BFGS and dBN for approximating function \cref{Example2eq}.}    \label{example3BFGSdBN}
\end{figure}

The first test problem is the function
\begin{equation}\label{Example2eq}
    u(x) = \sqrt{x}.
\end{equation}
as the target function for problem \cref{min}, with $r(x) = 1$. We aim to test the performance of dBN  for least-squares approximation problems. \cref{example3BFGSdBN} presents a comparison between dBN and L-BFGS. In this comparison, we utilized a Python L-BFGS implementation from `scipy.optimize'. The initial network parameters for the two algorithms were set to be the uniform mesh for $\bb^{(0)}$ with $\bc^{(0)}$ given by solving $\cref{lineaEqDF}$. The computation times for selected iterations are reported in \cref{tab:method_comparison} (CPU: 12th Gen Intel(R) Core(TM) i3-1215U, 1.20GHz).  In this example, our solver outperforms L-BFGS, achieving smaller losses in fewer iterations and less time.

\begin{figure}[htb!]
    \centering
    \subfigure[]
    {\includegraphics[width=0.45\textwidth]{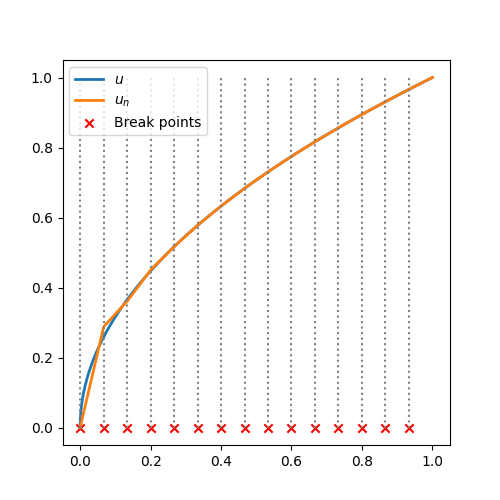}}
    \hspace{2em}
    \subfigure[]{
        \includegraphics[width=0.45\textwidth]{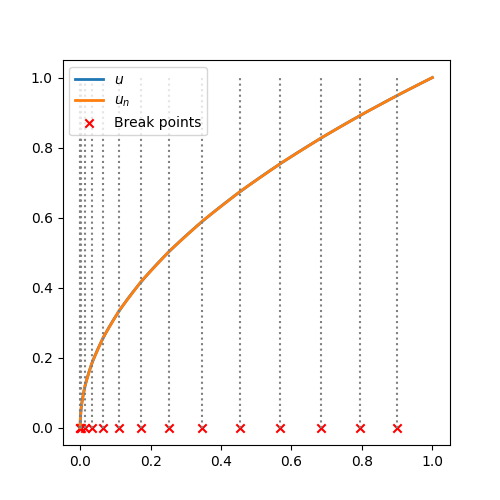}}
    \caption{(a) initial NN model with 15 uniform breakpoints, $J(u_n) = 5.64 \times 10^{-5}$, (b) optimized NN model with 15 breakpoints using dBN, 1000 iterations, $J(u_n) = 1.20 \times 10^{-7}$.} \label{example2DF}
\end{figure}

\begin{table}[ht]
    \centering
    \begin{tabular}{|c|c||c|c|c|c|}
        \hline
         \multicolumn{2}{|c||}{} & \multicolumn{4}{c|}{Computation time (seconds)} \\
        \hline
        Method& $n$ &  20 itr. & 40 itr. & 60 itr. & 80 itr. \\ 
        \hline
        dBN    & 24 & 0.16 & 0.29 & 0.43 & 0.56 \\
        L-BFGS & 24 & 0.66 & 1.39 & 2.10 & 2.80 \\
        dBN    & 48 & 0.39 & 0.75 & 1.09 & 1.50 \\
        L-BFGS & 48 & 2.67 & 5.01 & 7.59 & 10.11 \\
        \hline
    \end{tabular}
        \caption{Times at various iterations of dBN and L-BFGS for approximating function \cref{Example2eq}.}\label{tab:method_comparison}
\end{table}
\Cref{example2DF} (a) illustrates the neural network approximation of the function in \cref{Example2eq}, obtained using uniform breakpoints and determining the linear parameters through the solution of \cref{lineaEqDF}. Clearly, it is more optimal to concentrate more mesh points on the left side, where the curve is steeper. The dBN method is capable of making this adjustment, as illustrated in \cref{example2DF} (b). The loss functions confirm that the approximation improves substantially when the breakpoints are allocated according to the steepness of the function.

\subsection{Exponential Solution}\label{sectionex1} The second test problem involves the function
\begin{equation}\label{Example1eq}
u(x) = x \left( \exp \left( -\frac{{(x - \frac{1}{3})^2}}{{0.01}} \right) - \exp \left( -\frac{{4}}{{9 \times 0.01}} \right) \right),
\end{equation}
serving as a solution of \cref{pde} for $a(x) = r(x)=1$ and $\alpha = \beta = 0$.

\begin{figure}[htb!]
    \centering
    \subfigure[Relative error $e_n$ versus number of iterations with 24 breakpoins. Final relative errors: L-BFGS - 0.183, dBN - 0.076]
    {\includegraphics[width=0.45\textwidth]{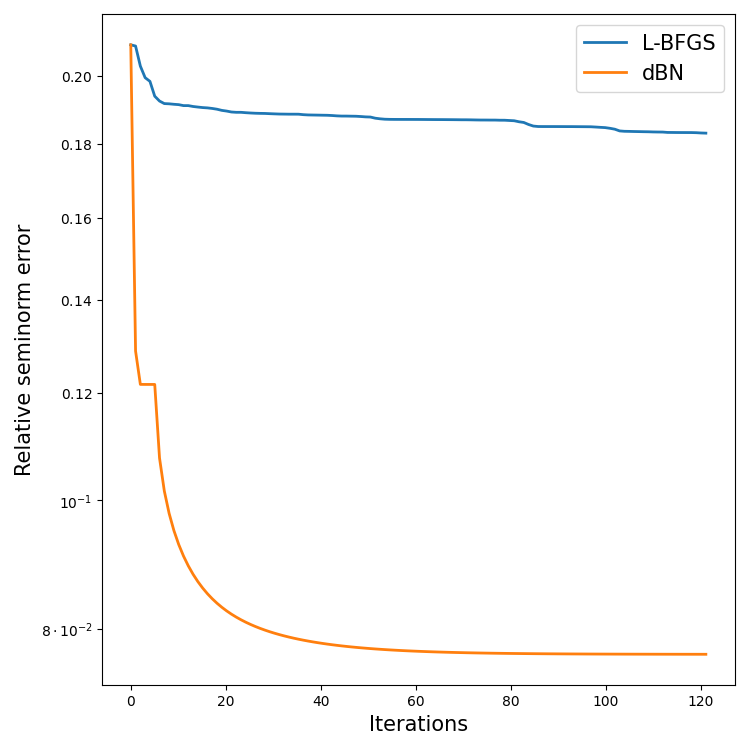}}
    \hspace{2em}
    \subfigure[Relative error $e_n$ versus number of iterations with 48 breakpoints. Final relative errors: L-BFGS - 0.099, dBN - 0.056]{
        \includegraphics[width=0.45\textwidth]{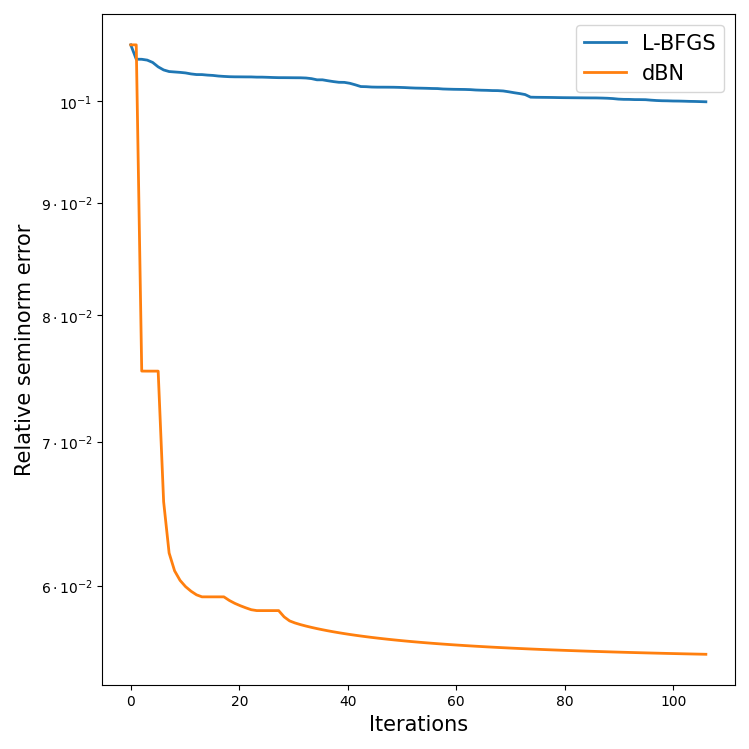}}
    \caption{Comparison between L-BFGS and dBN for approximating function \cref{Example1eq}.}\label{example1BFGSdB}
\end{figure}

Similarly to \cref{example3BFGSdBN}, we start by comparing our solver with L-BFGS. The initial network parameters for both algorithms were set to be the uniform mesh for $\bb^{(0)}$, with $\bc^{(0)}$ given by the exact solution of equation \cref{linear_eq}. \Cref{tab:method_comparison2} reports the times at various iterations. From \cref{example1BFGSdB} and \cref{tab:method_comparison2}, we observe that dBN outperforms L-BFGS in both accuracy and computation time.

\begin{table}[ht]
    \centering

    \begin{tabular}{|c|c||c|c|c|c|}
        \hline
         \multicolumn{2}{|c||}{} & \multicolumn{4}{c|}{Computation time (seconds)} \\
        \hline
        Method& $n$ &  20 itr. & 40 itr. & 60 itr. & 80 itr. \\ 
        \hline
        dBN    & 24 & 0.14 & 0.30 & 0.45 & 0.60 \\
        L-BFGS & 24 & 0.80 & 1.48 & 2.16 & 2.90 \\
        dBN    & 48 & 0.35 & 0.71 & 1.08 & 1.43 \\
        L-BFGS & 48 & 2.97 & 5.95 & 8.79 & 10.09 \\
        \hline
    \end{tabular}
        \caption{Times at various iterations of dBN and L-BFGS for approximating function \cref{Example1eq}.}\label{tab:method_comparison2}
\end{table}

\begin{figure}[htb!]
  \centering 
  \subfigure[]{ 
    \includegraphics[width=0.31\linewidth]{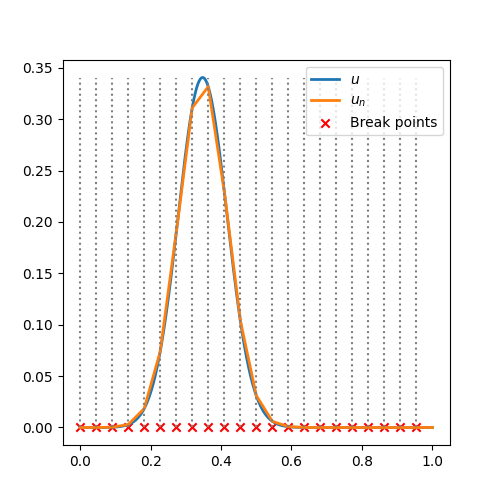}}
    \subfigure[]{ 
\includegraphics[width=0.31\linewidth]{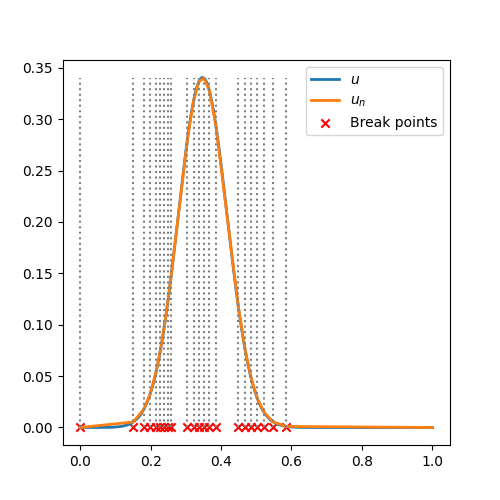}}
    \subfigure[]{\includegraphics[width=0.31\linewidth]{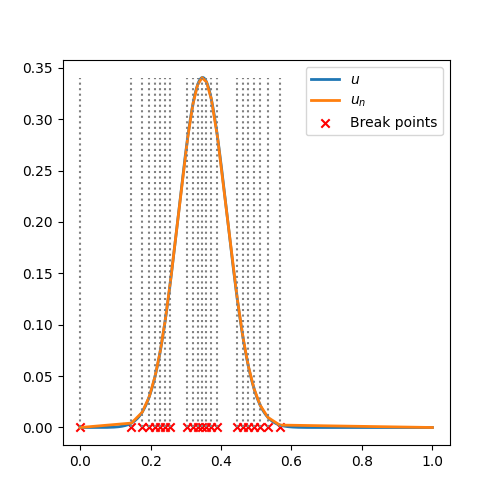}
    }
    \caption{(a) initial NN model with 22 uniform breakpoints, $e_n = 0.228$, (b) optimized NN model with 22 breakpoints, 500 iterations, $e_n = 0.092$, (c) adaptive approximation ($n = 8, 11, 14, 18, 22$), $e_n =  0.083$.}\label{ex1Figure}
\end{figure}

\cref{ex1Figure} (a) shows the initial neural network approximation of the function in \cref{Example1eq}, obtained by using uniform breakpoints and determining the linear parameters through the solution of \cref{linear_eq}. The approximation generated by dBN is shown in \cref{ex1Figure} (b), while \cref{ex1Figure} (c) illustrates the approximation obtained by employing dBN with adaptivity (AdBN). Notably, in both cases, the breakpoints are moved, and the approximation enhances the initial approximation.

Theoretically, from \cref{Err}, $\frac{1}{n}$ is the order of convergence of approximating a solution \cref{Example1eq} by functions in ${\cal M}_n$. However, since \cref{min_pde} is a non-convex optimization problem, the existence of local minimums makes it challenging to achieve this order. Therefore, given the neural network approximation $u_n$ to $u$ provided by the dBN method, assume that 
\begin{equation*}
    e_n = \left(\frac{1}{n}\right)^r,
\end{equation*}
for some $r > 0$. As in \cite{cai2024fast}, we can use the AdBN method to improve the order of convergence of the dBN method (achieve an $r$ closer to 1).

\cref{comparison_adapt} illustrates adaptive dBN (AdBN) starting with 18 neurons, refining 9 times, and reaching a final count of 54 neurons. The stopping tolerance was set to $\epsilon = 0.12$. The recorded data in \cref{comparison_adapt} includes the relative seminorm error $e_n$ and the relative error estimator $\xi_n$ for each iteration of the adaptive process. Additionally, \cref{comparison_adapt} provides the results for dBN with fixed 50 and 54 neurons. Comparing these results to the adaptive run with the same number of neurons, we observe a significant improvement in rate, error estimator, and seminorm error within the adaptive run.


\begin{table}[htb!]
    \begin{center}
        \begin{tabular}{ |p{3cm}||p{2cm}|p{1.5cm}|p{1.5cm}|}
        \hline
            NN ($n$ breakpoints)&  $e_n$&  $\xi_n$ & $r$\\
             \hline
            Adaptive (18) & $1.11 \times10^{-1}$ & 0.371& 0.777  \\ 
            Adaptive (22) & $8.43 \times10^{-2}$ & 0.335& 0.813\\ 
            Adaptive (26) & $6.99 \times10^{-2}$ & 0.263 & 0.827\\ 
            Adaptive (30) & $6.14 \times10^{-2}$& 0.231 & 0.828 \\ 
            Adaptive (34) & $5.55 \times 10^{-2}$ & 0.183 & 0.827\\
            Adaptive (38) & $4.84 \times 10^{-2}$ & 0.159 & 0.839\\ 
            Adaptive (42) & $4.35\times10^{-2}$ & 0.136 &0.844\\
            Adaptive (46) & $3.97\times10^{-2}$ & 0.130 &0.848\\
            Adaptive (50) & $3.70\times10^{-2}$ & 0.122 &0.847\\
            Adaptive (54) & $3.40\times10^{-2}$ & 0.114 &0.852\\
            \hline
            Fixed (50) & $4.74\times10^{-2}$ & 0.147 & 0.783  \\
            Fixed (54) & $4.24\times10^{-2}$ & 0.121& 0.796\\
             \hline
        \end{tabular}
        \caption{Comparison of an adaptive network with fixed networks for relative error $e_n$, relative error estimators $\xi_n$, and powers $r$.}
        \label{comparison_adapt}
    \end{center}
\end{table}

\subsection{Singularly Perturbed Reaction-Diffusion Equation}\label{sing_perturb}
The third test problem is a singularly perturbed reaction-diffusion equation:
\begin{equation}\label{pde2}
    \left\{\begin{array}{lr}
        -\varepsilon^2u^{\prime\prime}(x) +u(x)=f(x), & x\in I=(-1,1),\\[2mm]
        u(-1) = u(1) = 0.
    \end{array}\right.
\end{equation}
For $f(x) = -2\left(\varepsilon -4x^2\tanh{\left(\frac{1}{\varepsilon}(x^2 - \frac{1}{4})\right)}\right)\left(1/\cosh{\left( \frac{1}{\varepsilon}(x^2 - \frac{1}{4})\right)}\right)^2 + \tanh{\left(\frac{1}{\varepsilon}(x^2 - \frac{1}{4})\right)} - \tanh{\left(\frac{3}{4\varepsilon}\right)}$, problem \cref{pde2} has the following exact solution
\begin{equation}\label{ex2pde}
    u(x) = \tanh{\left(\frac{1}{\varepsilon}(x^2 - \frac{1}{4})\right)} - \tanh{\left(\frac{3}{4\varepsilon}\right)}.
\end{equation}

\begin{figure}[htb!]
    \centering
    \subfigure[]
    {\includegraphics[width=0.45\textwidth]{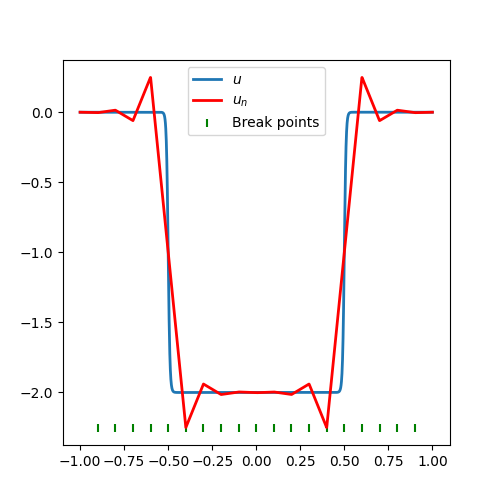}}
    \hspace{2em}
    \subfigure[]{
        \includegraphics[width=0.45\textwidth]{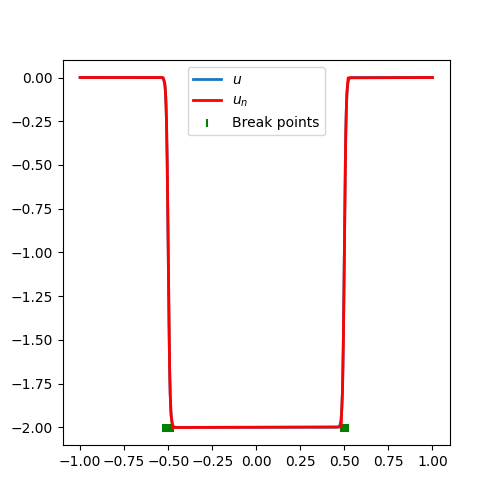}}
    \caption{For $\nu = \varepsilon^2 = 10^{-4}$: (a) initial NN model with 20 uniform breakpoints, $e_n = 0.935$, (b) optimized NN model with 20 breakpoints, 200 iterations, $e_n =0.130$.}\label{example2DR}

\end{figure}

For some $\nu = \varepsilon^2$, these problems exhibit interior layers that make them challenging for mesh-based methods such as finite element and finite difference, leading to overshooting and oscillations. For $\nu = 10^{-4}$, \cref{example2DR} illustrates the neural network approximation of the function described in \cref{ex2pde}, using uniform breakpoints (a) and employing dBN to adjust the breakpoints (b). An interesting observation is that the resulting approximation from dBN does not exhibit overshooting or oscillations. This confirms that dBN is capable of successfully adjusting the breakpoints and may have the potential to accurately approximate solutions with boundary and/or interior layers.

It is worth mentioning that the relative $L^2$-norm error of the approximation depicted in \cref{example2DR} (b) is $2.12\times 10^{-3}$. In \cite{CAI2020109707}, similar $L^2$ errors were obtained using deep neural networks with $2962$ parameters. In our case, the number of parameters is only $41$.

The resulting relative errors (measured in the $H^1$ seminorm as given in \cref{seminorm_error}) obtained after using dBN for 32 breakpoints and various values of $\nu$ are shown in \cref{comparison_epsilon}. For each value of $\nu$, dBN considerably improves the initial approximation, and the error does not vary significantly with different values of $\nu$.

\begin{table}[htb!]
    \begin{center}
        \begin{tabular}{ |p{1cm}||p{2cm}|p{2cm}|}
        \hline
            $\nu$& $e_n$ (initial)   &$e_n$ (dBN) \\
             \hline
            $10^{-2}$& $1.63\times 10^{-1} $&$6.72 \times10^{-2}$ \\ 
            $10^{-3}$& $5.53\times10^{-1}$&$8.08 \times10^{-2}$ \\ 
            $10^{-4}$& $8.89\times10^{-1}$ &$7.65 \times10^{-2}$ \\ 
            $10^{-5}$& $9.69\times10^{-1}$&$8.58 \times10^{-2}$ \\
            $10^{-6}$&  $9.90\times10^{-1}$&$8.92 \times10^{-2}$ \\

             \hline
        \end{tabular}
        \caption{Relative error $e_n$ 
        for $32$ breakpoints: initial approximation on uniform breakpoints
        and approximation after 500 iterations of dBN.}
        \label{comparison_epsilon}
    \end{center}
\end{table}

\begin{figure}[htb!]
    \centering
    \subfigure[]
    {\includegraphics[width=0.45\textwidth]{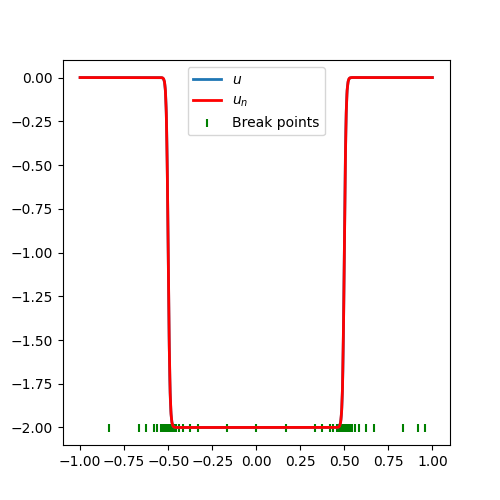}}
    \hspace{2em}
    \subfigure[]{
        \includegraphics[width=0.45\textwidth]{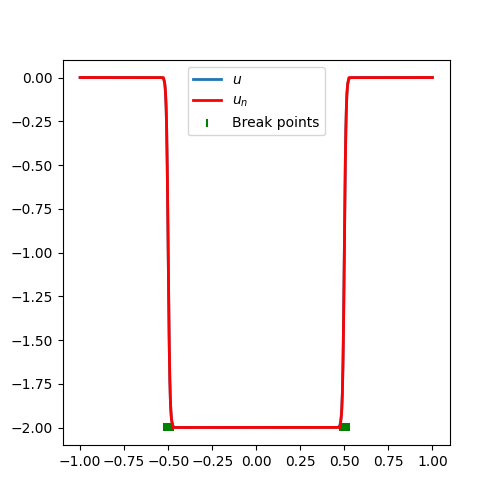}}
    \caption{For $\nu = \varepsilon^2 = 10^{-4}$: (a) aFEM approximation  ($n = 12, 13, 18, 26, 34, 42, 54$), $e_n = 0.133$, (b) AdBN approximation ($n = 12, 16, 17, 22, 28, 29, 30$), $e_n = 0.079$.}\label{example22DR}


\end{figure}

We also present the results using an adaptive finite element method (aFEM). \cref{comparison_adapt2} compares AdBN and aFEM. In this experiment, both AdBN and aFEM start with 12 breakpoints and are refined 6 times. The same marking strategy from \cref{marking} is used to determine the intervals for refinement in aFEM. After each refinement, the linear parameters were computed by solving equation \cref{linear_eq}.  At the third refinement of AdBN, with 22 breakpoints, we obtain a more accurate approximation than the one at the final refinement of aFEM with 54 breakpoints. The resulting approximations are shown in \cref{example22DR}.

\begin{table}[h!]
    \begin{center}
        \begin{tabular}{ |p{4cm}||p{2cm}|p{1.5cm}|}
        \hline
            Method ($n$ breakpoints)&  $e_n$\\
             \hline 
            AdBN (12) & $2.89 \times10^{-1}$\\
            AdBN (16) & $1.71 \times10^{-1}$\\
            AdBN (17) & $1.66 \times10^{-1}$\\
            AdBN (22) & $1.32 \times10^{-1}$ \\
            AdBN (28) & $9.42 \times10^{-2}$ \\
            AdBN (29) & $8.10 \times10^{-2}$\\
            AdBN (30) & $7.88 \times10^{-2}$\\        
            \hline
            aFEM (12) & $9.63 \times 10^{-1}$\\ 
            aFEM (13) & $9.63 \times 10^{-1}$\\ 
            aFEM (18) & $9.19 \times 10^{-1}$\\
            aFEM (26) & $8.18 \times 10^{-1}$\\
            aFEM (34) & $5.77 \times 10^{-1}$\\
            aFEM (42) & $2.73 \times 10^{-1}$\\
            aFEM (54) & $1.33 \times 10^{-1}$\\

             \hline
        \end{tabular}
        \caption{Comparison of AdBN and aFEM for relative errors $e_n$.  
        }
        \label{comparison_adapt2}
    \end{center}
\end{table}

\section{Discussion and Conclusion}The shallow Ritz method improves the order of approximation for one-dimensional non-smooth elliptic PDEs. The dBN method is an efficient iterative approach for solving the computationally intensive non-convex optimization problem arising from finding the optimal breakpoints. Extending this method to diffusion-reaction problems presents additional difficulties.

First, unlike the finite element mass matrix, the NN mass matrix is dense and very ill-conditioned. This difficulty is overcome in one dimension through a special factorization of the mass matrix, which was done using both algebraic and geometrical approaches. This factorization enables the ${\cal O}(n)$ computational cost for the inversion of the mass matrix.

Second, the optimality conditions for the non-linear parameters are no longer nearly decoupled non-linear algebraic systems, which leads to dense Hessian matrices. Additionally, the Hessians may be singular. To address this, we define a reduced non-linear system in which the Hessian matrix is invertible. The inverse of the reduced system’s Hessian matrix can be computed in ${\cal O}(n)$ operations by exploiting the fact that the inverse of the stiffness matrix is tri-diagonal. Consequently, the resulting damped block Newton (dBN) method is implemented with a computational cost of only ${\cal O}(n)$ per iteration.


Overall, the numerical results demonstrate the efficiency of our method in terms of not only the number of iterations but also the cost per iteration, making a compelling case to pursue the construction of similar solvers for higher dimensional problems. Of particular interest is the application of dBN methods to the singularly perturbed reaction-diffusion problem. For a fixed number of mesh points $n$, dBN appears to achieve an accuracy independent of the diffusion coefficient $\varepsilon^2$. Furthermore, when adding in adaptivity, AdBN seems to be comparable to FE methods using mesh refinement.

The ideas developed for the dBN method in diffusion-reaction problems are readily applicable to least-squares approximation problems. For both problems (elliptic PDEs and least-squares approximation), local convergence can be guaranteed; a detailed convergence analysis will be presented in a forthcoming paper.

\appendix
\section{Inverting the Mass Matrix}\label{appendixA}
\subsection{Algebraic Approach}\label{appendixA_algebrajc}
This section derives an inverse formula of the mass matrix through a decomposition into two matrices. 
The decomposition is based on the fact that matrices with the structure of ${\cal M}$ in \cref{specialMatrix} have tri-diagonal inverses.

For $1\leq i\leq j \leq n+1$, let $m_{ij}$ be the $(i,j)$-element of the mass matrix $M_r(\bb)$, then
\begin{eqnarray*}
    m_{ij} &=&  \int_0^1 r(x)\sigma(x-b_{i-1})\sigma(x-b_{j-1})dx = \int_{b_{j-1}}^1 r(x) (x-b_{i-1})(x-b_{j-1})dx \\
    &=& \int_{b_{j-1}}^1 r(x) \left(x-1\right)\left(x-b_{j-1}\right) dx + (1-b_{i-1})\int_{b_{j-1}}^1 r(x)\left(x - b_{j-1} \right)dx \equiv m^1_{ij} + m^2_{ij},
\end{eqnarray*}
which implies the following decomposition
\begin{equation*}\label{decomp}
M_r(\bb) = M^1(\bb) + M^2(\bb) \equiv \left(m_{ij}^1\right)_{n\times n} + \left(m_{ij}^2\right)_{n\times n}.
\end{equation*}
Both $M^1(\bb)$ and $M^2(\bb)$ have the same structure as ${\cal M}$ in \cref{specialMatrix} with
\[
m_{ij}^1=\beta^1_{\max\{i,j\}}  \quad\mbox{and}\quad m_{ij}^2= \alpha_{min\{i,j\}}^2\beta^2_{max\{i,j\}},
\]
where
\[
\beta^1_k=\int_{b_{k-1}}^1 r(x) \left(x-1\right)\left(x-b_{k-1}\right) dx, \quad
\alpha^2_k=1-b_{k-1}, \quad\mbox{and}\quad
\beta^2_k=\int_{b_{k-1}}^1 r(x)\left(x - b_{k-1} \right)dx.
\]
\begin{proposition}\label{prop_invM}
The inverse of the mass matrix $M_r(\bb)$ is given by
\begin{equation}\label{inv-M}
M_r(\bb)^{-1}=M^2(\bb)^{-1}(M^2(\bb)^{-1} + M^1(\bb)^{-1})^{-1}M^1(\bb)^{-1}.
\end{equation}
\end{proposition}
\begin{proof}
\cref{inv-M} is a direct consequence of the fact that
\[
M_r(\bb) = M^1(\bb) \left(M^2(\bb)^{-1} + M^1(\bb)^{-1} \right) M^2(\bb).
\]
\end{proof}

\begin{remark}
Since $M^1(\bb)^{-1}$ and $M^2(\bb)^{-1}$ are tri-diagonal, so is 
$M^1(\bb)^{-1}\! + \!M^2(\bb)^{-1}$. Hence, $M_r(\bb)^{-1}$ in \cref{inv-M} applied to any vector can be computed in ${\cal O}(n)$ operations.
\end{remark}

\subsection{Geometric Approach}\label{appendixA_geometric}This section presents another way to invert the mass matrix, based on a factorization of $M_r(\bb)$ into the product of three  matrices easy to invert. The factorization arises from expressing the global ReLU basis functions in terms of local discontinuous basis functions. 

To this end, for $k = 0, \ldots, n$, let $I_k = [b_k, b_{k+1})$ and  define the local basis functions 
\begin{align*} \label{eqn-disc-phi}
  \varphi_k^0(x) = \left\{\begin{array}{cl}
  1, & x \in I_k, \\ 0, & \text{otherwise} \end{array}\right. \quad\mbox{and}\quad
  \varphi_k^1 (x) =
  \left\{\begin{array}{cl}
  h_k^{-1}(x-b_k), & x \in I_k, \\ 0, & \text{otherwise} \end{array}\right. .
\end{align*}
Since $\sum\limits_{k=0}^n\varphi_k^0(x)\equiv  1$ in $I$, we have
\begin{equation}\label{inclusion}
\mbox{span}\left\{1, \sigma(x-b_0),\ldots, \sigma(x-b_n) \right\} \subset \mbox{span} \left\{\varphi^0_k(x) \right\}_{k=0}^{n} \bigcup \mbox{span} \left\{\varphi^1_k(x) \right\}_{k=0}^{n}.
\end{equation}
Set 
\begin{equation}\label{psi}
\bpsi(x)=(\psi_0(x), \ldots,\psi_n(x))^T \quad\mbox{and}\quad \bphi_{i}(x) = (\varphi_0^i(x), \ldots, \varphi^i_n(x))^T,
\end{equation}
where $\psi_k(x)=\sigma(x-b_k)$; and let $\bD(\mathbf{h}) = \mbox{diag}(h_0,\ldots, h_{n})$, 
\begin{align*}\label{matricesG&D}
  G = \left(
  \begin{array}{cccc}
      ~1 &      &    & \\
    \!-1 &   ~1 &    & \\
         &   \ddots   &     \ddots & \\
         &      &     \!-1 &~1
  \end{array}
  \right)_{\!\!{n+1}\times {n+1}},
  \quad\mbox{and}\quad
  G^{-1} = \left(
  \begin{array}{cccc}
      1 &      &    & \\
    1 &   ~1 &    & \\
     \vdots    &  \vdots    &     \ddots & \\
    1     & 1     & \ldots    &~1
  \end{array}
  \right)_{\!\!{n+1}\times {n+1}}.
\end{align*}

\begin{lemma}\label{basisFunctions}
There exist mappings $B_0: \R^{n+1} \to \R^{n+1}$ and $B_1: \R^{n+1} \to \R^{n+1}$ such that 
    \begin{equation}\label{trans}
        \bpsi = B_0\bphi_0+B_1\bphi_1.
    \end{equation} 
Moreover, we have 
\[
B_0=G^{-T}\bD(\mathbf{h}) \left(G^{-T} -I\right) \quad\mbox{and}\quad B_1 =G^{-T}\bD(\mathbf{h}),
\]
where $I$ is the $(n+1)$-order identity matrix.
\end{lemma}
\begin{proof}
\cref{inclusion} implies that there exist $B_0$ and $B_1$ such that \cref{trans} is valid. To determine $B_0$ and $B_1$, for any $\bc=(c_0,\ldots,c_n)^T\in \R^{n+1}$, let $v(x)=\bc^T\bpsi(x)$, then
\[
v(x)=\bc^T\bpsi(x)=\bc^T B_0\bphi_0(x) + \bc^TB_1\bphi_1(x).
\]
On each $I_k$, using the fact that $v^\prime(x)$ and $\bc^T B_0\bphi_0(x)$ are both constant, we have
\[
  \left(\bc^TB_1\right)_{k+1}=v(b_{k+1}) - v(b_{k}) 
  = \sum_{i=0}^k c_i \left(\sum_{j=i}^k h_j\right) -
     \sum_{i=0}^{k-1} c_i \left(\sum_{j=i}^{k-1} h_j\right) 
  = \sum_{i=0}^k c_i h_k =  \left(\bD(\bh)G^{-1} \mathbf{c}\right)_{k+1}
\]
which, together with the arbitrariness of $\bc$, implies that $B_1=G^{-T}\bD(\bh)$.

By the definitions of $\bphi_0(x)$ and $\bphi_1(x)$ and the fact that $v(b_0)=0$, we have
\begin{align*}
  \left(\bc^TB_0\right)_{k+1} &= v(b_{k}) 
  = \left(v(b_{k}) - v(b_{k-1})\right) + \left(v(b_{k-1}) - v(b_{k-2})\right) + \cdots + \left(v(b_{1}) - v(b_{0})\right) \\
  &= \left(\bD(\bh)G^{-1} \mathbf{c}\right)_{k} + \left(\bD(\bh)G^{-1} \mathbf{c}\right)_{k-1} + \cdots +\left(\bD(\bh)G^{-1} \mathbf{c}\right)_1 =  \left( (G^{-1} - I)\bD(\mathbf{h}) G^{-1} \bc \right)_k,
\end{align*}
which, together with the arbitrariness of $\bc$, implies that $B_0=G^{-T}\bD(\mathbf{h}) \left(G^{-T} -I\right)$. This completes the proof of the lemma.
\end{proof}

For $i, j = 0, 1$, let 
\[
\bD_{ij}(r) =  \displaystyle\int_0^1 r(x)\bphi_i \bphi_j^T dx .
\]
For $k=0,1,2$, let 
\begin{equation}\label{Dr}
\bD_r(\bs^k)=\text{diag}(s_0^k(r),\dots,s_n^k(r)) \quad\mbox{with}\quad s_i^k(r)=\int_{b_i}^{b_{i+1}}r(x)(x-b_i)^k\,dx. 
\end{equation}
Then, together with $\bD(\bh) = \mbox{diag}(h_0,\ldots, h_{n})$, it is easy to see that
\begin{equation*}\label{Dij}
\bD_{00}(r)=\bD_r(\bs^0), \quad \bD_{01}(r)=\bD_{10}(r)=\bD(\bh)^{-1}\bD_r(\bs^1), \quad \mbox{ and }\,\, \bD_{11}(r)=\bD(\bh)^{-2}\bD_r(\bs^2) .
\end{equation*}

\begin{theorem}\label{thm: mass_factorization}
Let $Q=G\bD(\bh)^{-1}G$ and let
\[
T_{M_r}= (I-G^{T})\bD_{00}(r)(I-G) +(I-G^{T})\bD_{01}(r)G +G^T \bD_{10}(r)(I-G) + G^{T}\bD_{11}(r)G,
\]
then the mass matrix $M_r (\bb)$ defined in \cref{massMatrix} has the following factorization
    \begin{equation}\label{fac}
        M_r (\bb)= Q^{-T}T_{M_r}Q^{-1}.
    \end{equation}
\end{theorem}

\begin{proof}
By \cref{trans} and the fact that $B_0=B_1 (G^{-T}-I)$, we have
    \begin{align*}
    M_r(\bb) &= \int_0^1 r(x) \bpsi\bpsi^T\,dx = B_0 \bD_{00} B_0^T + B_0 \bD_{01} B_1^T + B_1 \bD_{10} B_0^T + B_1 \bD_{11} B_1^T \\
    & = B_1 \left\{ (G^{-T}-I)\bD_{00}(G^{-1}-I) + (G^{-T}-I)\bD_{01} + \bD_{10}(G^{-1}-I) + \bD_{11} \right\} B_1^T\\
    & = B_1G^{-T} \left\{ (I-G^{T})\bD_{00}(I-G) +(I-G^{T})\bD_{01}G +G^T \bD_{10}(I-G) + G^{T}\bD_{11}G \right\}
     G^{-1}B_1^T,
\end{align*}
which, together with the fact that $Q^{-1}=G^{-1}B_1^T$, implies \cref{fac}.
\end{proof}


\begin{remark}\label{rmk: coeff_factorization}
The transformation in \cref{trans} leads to a similar factorization of the stiffness matrix as  $A_a(\bb) = Q^{-T}T_{A_a}Q^{-1}$
   with $\displaystyle T_{A_a}= G^T \bD(\bh)^{-2}\bD_a(\bs^0)G$ being tri-diagonal, where $\bD_a(\bs^0)$ is defined similarly as in \cref{Dr}. Therefore, the sum of the stiffness matrix and the mass matrix satisfies that $A_a(\bb) + M_r(\bb) =  Q^{-T}(T_{A_a} + T_{M_r})Q^{-1}$. 
\end{remark}

\bibliographystyle{plain}
\bibliography{ref}

\clearpage

\end{document}